\documentclass[reqno]{amsart}
\usepackage{amsxtra}
\usepackage{amssymb,amsmath,amsfonts,latexsym,}
\usepackage{color}
\usepackage[all]{xy}
\usepackage{pb-diagram}
\numberwithin{equation}{section}

\theoremstyle{plain}

\newcommand{\cleqn}{\setcounter{equation}{0}}
\newcommand{\clth}{\setcounter{theorem}{0}}
\newcommand {\sectionnew}[1]{\section{#1}\cleqn\clth}

\newtheorem{theorem}{Theorem}[section]
\newtheorem{lemma}[theorem]{Lemma}
\newtheorem{proposition}[theorem]{Proposition}
\newtheorem{corollary}[theorem]{Corollary}
\newtheorem{definition}[theorem]{Definition}
\newtheorem{example}[theorem]{Example}
\newtheorem{remark}[theorem]{Remark}

\newcommand \lb[1]{\label{#1}}

\newcommand{\rmap}       {\longrightarrow}

\newcommand{\Ker}        {{\mathrm {ker}}}

\newcommand{\source}        {\mathsf{s}}
\newcommand{\target}         {{\mathsf{t}}}

\newcommand{\Cour}[1]      {[\![#1]\!]}

\newcommand{\Lie}        {\mathcal L}

\newcommand{\R}      {{\mathbb R}}

\begin{document}
\title{$B$-field transformations of Poisson groupoids}

\author{Cristi\'an Ortiz}
\address{Departamento de Matem\'atica\\
Universidade Federal do Paran\'a\\
Setor de Ciencias Exatas - Centro Politecnico
81531-990 Curitiba - Brazil.}
\email{cristian.ortiz@ufpr.br}

\date{}

\begin{abstract}
This is a survey devoted to the study of $B$-field transformations of multiplicative Poisson bivectors on a Lie groupoid $G$. We are concerned with $B$-fields given by  multiplicative closed $2$-forms on $G$. We extend the results in \cite{B} by viewing Poisson groupoids and their $B$-field symmetries as special instances of multiplicative Dirac structures \cite{Ortiz-Thesis}. We also describe such symmetries infinitesimally.
\end{abstract}

\maketitle

\hspace{4.4cm}\emph{Dedicado a la memoria de Leonor Gonz\'alez.}

\tableofcontents

\sectionnew{Introduction}\lb{intro}

The concept of Poisson groupoid was introduced by Alan Weinstein \cite{weinsteincoisotropic} as a common generalization of Poisson Lie groups \cite{D} and symplectic groupoids \cite{CDW}. A Poisson groupoid is a pair $(G,\pi_G)$ where $G$ is a Lie groupoid and $\pi_G$ is a Poisson bivector on $G$ suitably compatible with the groupoid multiplication. Poisson group(oid)s are interesting by themselves, but one of the main motivations for studying such geometric objects is due to their interaction with the theory of  integrable systems. On one hand, Poisson Lie groups help to study Hamiltonian properties of the group of dressing transformations of certain integrable systems where this group plays the role of hidden symmetries \cite{STS}, on the other hand Poisson group\emph{oids} are useful for describing the geometry of solutions of the classical dynamical Yang-Baxter equation \cite{EV}.

It was proven in \cite{Mac-Xu} that the compatibility of a Poisson bivector with a groupoid multiplication on $G$ is reflected infinitesimally by the fact that $(AG,A^*G)$ is a Lie bialgebroid. Here $AG$ is the Lie algebroid of $G$ and $A^*G$ is the dual vector bundle of $AG$.  
Hence, Lie bialgebroids are regarded as the \emph{infinitesimal data} of Poisson groupoids. The integration of Lie bialgebroids to Poisson groupoids was carried out in \cite{Mac-Xu2}. More precisely if $A$ is the Lie algebroid of a source simply connected Lie groupoid $G(A)$ and the pair $(A,A^*)$ is a Lie bialgebroid, then there exists a unique Poisson bivector $\pi_{G(A)}$ on $G(A)$ making the pair $(G(A),\pi_{G(A)})$ into a Poisson groupoid whose Lie bialgebroid is $(A,A^*)$. The relation between Poisson groupoids and Lie bialgebroids generalizes the well known Drinfeld's correspondence between Poisson Lie groups and Lie bialgebras \cite{D}.

Since Poisson groupoids are important examples of Poisson manifolds, it is interesting to study symmetries of these structures. In \cite{SW} it was observed that a Poisson structure $\pi$ on a smooth manifold $M$ can be transformed by a closed $2$-form $B\in\Omega^2(M)$ into another Poisson manifold as follows: look at the symplectic foliation determined by $\pi$ and add to the leafwise symplectic form the pull-back of $B$ to the symplectic leaves; if the resulting leafwise $2$-form is symplectic, then it is given by a Poisson structure $\pi_B$ on the same manifold $M$. The Poisson structure $\pi_B$ is called a \textbf{gauge transformation} of $\pi$ by the closed $2$-form $B\in\Omega^2(M)$ or a \textbf{$B$-field transformation} of $(M,\pi)$.  If the resulting $2$-form is not symplectic, then it does not come from a Poisson structure. Instead, the underlying geometry is the one determined by a Dirac structure \cite{courant,SW}. Dirac structures were introduced in \cite{courant} as a unified approach to Poisson structures, pre-symplectic forms and regular foliations. 

In this work we study $B$-field transformations of Poisson groupoids. In order to preserve the compatibility of the geometric structure with the groupoid multiplication, we restrict our attention to multiplicative $B$-fields, that is, $B$-field transformations given by multiplicative closed $2$-forms. Given a Lie groupoid $G$ over $M$ with source and target maps $\source,\target:G\rmap M$ and multiplication $m_G:G_{(2)}\rmap G$, where $G_{(2)}:=\{(g,h)\in G\times G\mid \target(h)=\source(g)\}$, we say that a $2$-form $\omega_G\in\Omega^2(G)$ is \textbf{multiplicative} if  

$$m^*_G\omega_G=\mathrm{pr}^*_1\omega_G+\mathrm{pr}^*_2\omega_G,$$

\noindent where $\mathrm{pr}_1,\mathrm{pr}_2:G_{(2)}\rmap G$ are the canonical projections. Our main goal is to understand how the infinitesimal data of a Poisson groupoid $(G,\pi_G)$ changes after applying a multiplicative $B$-field. This problem is closely related to the concept of IM-$2$-form on a Lie algebroid \cite{BCWZ}. Given a Lie algebroid $A$ over $M$ with anchor $\rho_A:A\rmap TM$,  a bundle map $\sigma:A\rmap T^*M$ is called an \textbf{IM-$2$-form} on $A$ if the following identities hold
\begin{align*}
\langle \sigma(a), \rho_{A}(b)\rangle=&-\langle \sigma(b),\rho_{A}(a)\rangle\\
\vspace{.2cm}
\sigma[a,b]_{A}=&\Lie_{\rho_{A}(a)}\sigma(b)-\Lie_{\rho_{A}(a)}\sigma(b)+d\langle \sigma(a),\rho_{A}(b)\rangle.
\end{align*}

\noindent for every $a,b\in\Gamma(A)$. Alternatively, a bundle map $\sigma:A\rmap T^*M$ is an IM-$2$-form on the Lie algebroid $A$ if and only if the $2$-form $\omega_A:=-\sigma^*\omega_{can}\in\Omega^2(A)$ induces a Lie algebroid morphism

$$\omega^{\sharp}_A:TA\rmap T^*A \qquad X\mapsto \omega_A(X,\cdot),$$

\noindent between the tangent and cotangent Lie algebroids. Here $\omega_{can}$ is the canonical symplectic form on $T^*M$ (see \cite{BCO}). Throughout this paper, such $2$-forms on a Lie algebroid are called \textbf{morphic}. As shown in \cite{BCWZ} multiplicative closed $2$-forms on a source simply connected Lie groupoid are in one-to-one correspondence with IM-$2$-forms on the Lie algebroid of $G$. Equivalently \cite{BCO}, multiplicative closed $2$-forms on $G$ are in bijection with morphic $2$-forms on the Lie algebroid of $G$. This will be useful for describing multiplicative $B$-field transformations at the infinitesimal level.

First we study $B$-field transformations of symplectic groupoids. We observe that a $B$-field transformation of a Poisson manifold induces naturally an IM-$2$-form \cite{BCWZ, BCO} on the Lie algebroid of the Poisson manifold. This allows us to use the techniques developed in \cite{BCO} to describe the symplectic groupoid integrating gauge transformations of Poisson manifolds, providing an alternative proof of the results in \cite{BRad}. We also observe that multiplicative $B$-field transformations of Poisson groupoids can be studied by using multiplicative Dirac structures \cite{Ortiz-Thesis}. These are geometric structures that unify both multiplicative Poisson bivectors and multiplicative closed $2$-forms. Concretely, a Dirac structure $L_G$ on a Lie groupoid $G$ is called \textbf{multiplicative} if $L_G\subseteq TG\oplus T^*G$ is a Lie subgroupoid of the direct sum Lie groupoid. It was proved in \cite{Ortiz-Thesis} that, for every Lie groupoid $G$ with Lie algebroid $AG$, a multiplicative Dirac structure $L_G$ on $G$ corresponds to a Dirac structure $L_{AG}$ on $AG$ suitably compatible with both the linear and algebroid structures on $AG$. Using this terminology, we prove that multiplicative closed $2$-forms on a Lie groupoid $G$ act on the space of multiplicative Dirac structures, via $B$-field transformations. We show that this action is infinitesimally described by the $B$-field action of the corresponding morphic $2$-form on the space of Dirac structures on $AG$ compatible with both the vector bundle and Lie algebroid structures on $AG$. We also show that if $(G,\pi_G)$ is a Poisson groupoid with Lie bialgebroid $(AG,A^*G)$, then a multiplicative $B$-field transformation of $(G,\pi_G)$ corresponds infinitesimally to a morphic $B$-field transformation of the linear Poisson structure $\pi_{AG}$ induced by the dual Lie algebroid $A^*G$. If the resulting $B$-field transformation gives rise to a Poisson structure, we conclude that this Poisson structure is also multiplicative and our infinitesimal description recovers the results proved in \cite{B}. One can also study Poisson groupoids equipped with a Poisson action of a Lie group. If we assume that this action has the additional property of being given by groupoid automorphisms, then we can describe these actions at the infinitesimal level. This will be done in a separate paper.

This work is organized as follows: In section 2 we review the main definitions and examples of Lie groupoids and Lie algebroids. We also recall the definition of the tangent and cotangent groupoids as well as the tangent and cotangent Lie algebroids. In section 3 we study both multiplicative Poisson bivectors and closed $2$-forms, explaining what is the infinitesimal data of multiplicative Poisson bivectors \cite{Mac-Xu, Mac-Xu2} and also saying how closed multiplicative $2$-forms are described infinitesimally according to \cite{BCWZ, BCO}. Then in section 4 we study symmetries of Poisson structures viewed as special instances of Dirac structures. In section 5 we 
introduce multiplicative Dirac structures, describing such structures infinitesimally. We use this description to provide a conceptually clear setting for studying $B$-field transformations of Poisson groupoids.

\subsection{Acknowledgments} This article is an expanded version of a talk given at  the Second Latin Congress on Symmetries in Geometry and Physics, held at Universidade Federal do Paran\'a, Curitiba, in December 13-17, 2010. The author thanks Eduardo Hoefel and Elizabeth Gasparin  for the invitation to this greatly enjoyable conference. I also thank Henrique Bursztyn and Madeleine Jotz for useful comments and suggestions about the manuscript. The author also thanks IMPA, Rio de Janeiro and ESI, Vienna, where part of this work was carried out.

\subsection{Notations}

Let $M$ be a smooth manifold. The tangent bundle of $M$ is denoted by $p_M:TM\rmap M$. We use $c_M:T^*M\rmap M$ to indicate the cotangent bundle of a smooth manifold. Vector bundles are always denoted by $q_A:A\rmap M$ and the fiber over a point $x\in M$ is $A_x$.


\section{Lie theory}
\subsection{Lie groupoids}

A \textbf{groupoid} is a small category where all morphisms are invertible. More concretely, a groupoid is determined by a set $M$ of objects, a set $G$ of morphisms and structural maps $\source,\target:G\rmap M$ called \textbf{source} and \textbf{target}, respectively; a unit section $\epsilon_M:M\rmap G$ an inversion map $i_G:G\rmap G$ and a partially defined multiplication $m_G:G_{(2)}\rmap G$, where $G_{(2)}=\{(g,h)\in G\times G\mid \target(h)=\source(g)\}$ is the set of composable groupoid pairs. We require that all these maps satisfy the axioms of a category. See \cite{Mac-book} for more details. In this case we say that $G$ is a groupoid over $M$.

A groupoid $G$ over $M$ is called a \textbf{Lie groupoid} if both $G$ and $M$ are smooth manifolds, the source and target maps $\source,\target:G\rmap M$ are surjective submersions, and all the other structural maps are smooth. The submersion condition on $\source,\target$ ensures that $G_{(2)}$ inherits the structure of smooth manifold, so the smoothness of $m_G:G_{(2)}\rmap G$ makes sense. Throughout this work we only consider Lie groupoids.



\begin{example}

A Lie group can be thought of a as Lie groupoid whose space of units consist of only one point.

\end{example}

\begin{example}\label{transformationgroupoid}
Let $H$ be a Lie group acting on a smooth manifold $M$. We endow $H\times M$ with a Lie groupoid structure over $M$ as follows. The source and target maps are defined by

$$\source(h,x)=x, \quad \target(h,x)=hx.$$

\noindent The multiplication is defined by $(h,h'x)(h',x)=(hh',x)$. The unit section is $\epsilon(x)=(e,x)$ where $e\in H$ is the identity element. Finally the inversion map is defined by $i(h,x)=(h^{-1},hx)$. These maps define a Lie groupoid structure on $H\times M$, called the \textbf{transformation groupoid}. We usually denote the transformation groupoid by $H\ltimes M$. See \cite{Mac-book} for more details.
 
\end{example}

The previous examples show that Lie groupoids generalize not only Lie groups, but also Lie group actions.




Let $G_1$ and $G_2$ be Lie groupoids over $M_1$ and $M_2$, respectively. A \textbf{morphism} of Lie groupoids is a pair $(\Phi,\varphi)$ of smooth maps $\Phi:G_1\rmap G_2$ and $\varphi:M_1\rmap M_2$, compatible with all the structural maps, in the sense that $\Phi$ is a functor between the underlying categories. When $\Phi$ is injective we say that $G_1$ is a \textbf{Lie subgroupoid} of $G_2$.

\subsection{Lie algebroids}

A \textbf{Lie algebroid} is a vector bundle $q_A:A\rmap M$ equipped with a Lie bracket $[\cdot,\cdot]_A$ on the space of smooth sections $\Gamma_M(A)$ and a vector bundle map $\rho_A:A\rmap TM$ called the \textbf{anchor}, such that  

$$[a,fb]_A=f[a,b]_A+(\Lie_{\rho_A(a)}f)b,$$

\noindent for every $a,b\in\Gamma_M(A)$ and $f\in C^{\infty}(M)$.

\begin{example}
If $A$ is a Lie algebroid over a point we recover the notion of Lie algebra.

\end{example}

\begin{example}\label{transformationalgebroid}
Let $\mathfrak{h}$ be a Lie algebra acting on a smooth manifold $M$. That is, there exists a Lie algebra morphism

\begin{align*}
 \mathfrak{h}&\rmap\mathfrak{X}(M)\\
a&\mapsto a_M.
\end{align*}
 
\noindent We endow the trivial bundle $A_{\mathfrak{h}}=\mathfrak{h}\times M$ with the structure of a Lie algebroid over $M$. The anchor map is defined by 

\begin{align*}
 \rho:\mathfrak{h}\times M&\rmap TM\\
(a,x)&\mapsto a_{M}(x).
\end{align*}

\noindent The Lie bracket $[\cdot,\cdot]_{A_{\mathfrak{h}}}$ on $\Gamma(A_{\mathfrak{h}})\cong C^{\infty}(M)\bigotimes \mathfrak{h}$ is given by

$$[a,b]_{A_{\mathfrak{h}}}:=[a,b],$$

\noindent for $a,b\in\mathfrak{h}$, and we extend it by requiring the Leibniz rule. The bundle $A_{\mathfrak{h}}\rmap M$ with this Lie algebroid structure is referred to as the \textbf{transformation Lie algebroid}. See \cite{Mac-book} for more details.

\end{example}

Just as in the groupoid case, these examples show that Lie algebroids generalize both Lie algebras and their actions.




A Lie algebroid defines a differential graded algebra $\Omega^{\bullet}(A):=\Gamma_M(\wedge^{\bullet}A^*)$ with a degree $1$ operator $d_A:\Omega^{\bullet}(A)\rmap \Omega^{\bullet+1}(A)$ given by

\begin{align*}\label{algebroiddifferential}
d_A(\xi)(a_1,...,a_{p+1})=&\sum_{i=0}^{p}(-1)^{p+1}\Lie_{\rho_A(a_i)}\xi(a_1,...,\hat{a}_i,...a_{p+1})+\\ &+\sum_{i<j}(-1)^{i+j}\xi([a_i,a_j]_A,a_1,...,\hat{a}_i,...,\hat{a}_j,...,a_{p+1}).
\end{align*}

\noindent A direct computation shows that $d^2_A=0$ and that $d_A$ is a derivation of the exterior product. The operator $d_A$ determines the algebroid structure on $A$ in the sense that there is a one-to-one correspondence between degree $1$ derivations $D$ on $\Omega^{\bullet}(A)$ such that $D^2=0$ and Lie algebroid structures on the vector bundle $q_A:A\rmap M$. See for instance \cite{dufour}.

The latter description of Lie algebroids is useful to define morphisms between Lie algebroids. Let $A_1$ and $A_2$ be Lie algebroids over $M_1$ and $M_2$, respectively. A bundle map $\phi:A_1\rmap A_2$ covering $\varphi:M_1\rmap M_2$ induces an operator $\tilde{\phi}^*:\Gamma(A^*_2)\rmap \Gamma(A^*_1)$ defined by

$$(\tilde{\phi}^*\xi)_{x}a_x:=(\xi\circ \varphi)(x)\phi(a_x),$$

\noindent where $\xi\in \Gamma(A^*_2)$, $a\in\Gamma(A_1)$ and $x\in M$. This operator comes from the base preserving bundle map $\phi^*:\varphi^*(A^*_2)\rmap A^*_1$. We can extend this operator to higher wedge powers yielding an operator $\tilde{\phi}^*:\Omega^{\bullet}(A_2)\rmap \Omega^{\bullet}(A_1)$. We say that the bundle map $\phi:A_1\rmap A_2$ is a Lie algebroid \textbf{morphism} if the induced operator $\tilde{\phi}^*:\Omega^{\bullet}(A_2)\rmap \Omega^{\bullet}(A_1)$ is a chain map.

If $\phi:A_1\rmap A_2$ is injective, we say that $A_1$ is a \textbf{Lie subalgebroid} of $A_2$.

\begin{example}(Lie algebroids vs Poisson geometry)\label{poissonvsalgebroid}

Recall that a Poisson structure on a smooth manifold $M$ consists on a Lie bracket 

$$\{\cdot,\cdot\}:C^{\infty}(M)\times C^{\infty}(M)\rmap C^{\infty}(M),$$ 

\noindent such that $\{f,\cdot\}:C^{\infty}(M)\rmap C^{\infty}(M)$ is a derivation for every $f\in C^{\infty}(M)$. Equivalently, a Poisson structure on $M$ is given by a smooth \textbf{bivector} on $M$ (i.e. $\pi$ is a smooth section of the exterior bundle $\wedge^2 TM\rmap M$) such that $[\pi,\pi]=0$, where $[\cdot,\cdot]$ denotes the Schouten bracket of multivector fields \cite{dufour}. The correspondence is given by $\{f,g\}=\pi(df,dg)$. The pair $(M,\pi)$ is referred to as a \textbf{Poisson manifold}. Alternatively, a Poisson structure on $M$ corresponds to a Lie algebroid structure on the cotangent bundle $c_M:T^*M\rmap M$, where the anchor map is $\pi^{\sharp}:T^*M\rmap TM$, defined by

\begin{equation}\label{pisharp}
\beta(\pi^{\sharp}(\alpha)):=\pi(\alpha,\beta),
\end{equation}

\noindent where $\alpha,\beta\in T^*M$. The Lie bracket on $\Gamma(T^*M)=\Omega^1(M)$ is defined by

\begin{equation}\label{bracketT*M}
[\alpha,\beta]_{\pi}=\Lie_{\pi^{\sharp}(\alpha)}\beta - \Lie_{\pi^{\sharp}(\beta)}\alpha - \mathrm{d}\pi(\alpha,\beta).
\end{equation}

\noindent Thus, Lie algebroids arise naturally in Poisson geometry. Throughout this work, the Lie algebroid determined by a Poisson structure $\pi$ on $M$ will be denoted by $(T^*M)_{\pi}$.

\end{example}

It will be useful to give an alternative description of Lie algebroids in terms of Poisson structures. For that we observe that every Lie algebroid $A$ induces a Poisson structure on its dual bundle $A^*$ which is \textit{linear} in the sense that the space of fiberwise linear functions $C^{\infty}_{lin}(A^*)\cong\Gamma(A)\subseteq C^{\infty}(A^*)$ is a Poisson subalgebra. More precisely, the Poisson bracket $\{\cdot,\cdot\}_{A^*}$ on $C^{\infty}(A^*)$ is given by

\begin{itemize}

\item[i)] $\{f\circ q_{A^*}, g\circ q_{A^*}\}_{A^*}=0,$ for every $f,g\in C^{\infty}(M).$

\item[ii)] $\{a, q_{A^*}\circ f\}_{A^*}= q_{A^*}\circ (\Lie_{\rho_A(a)}f)$, for every $f\in C^{\infty}(M)$ and $a\in\Gamma(A)\cong C^{\infty}_{lin}(A^*).$

\item[iii)] $\{a,b\}_{A^*}=[a,b]_A$, for every $a,b\in \Gamma(A)\cong C^{\infty}_{lin}(A^*).$

\end{itemize}

It can be easily verified that this Poisson structure on $A^*$ determines completely the Lie algebroid structure on $A$. See e.g. \cite{CW}.

\subsection{The Lie functor}\label{liefunctor}

Here we recall the construction of the Lie functor from the category of Lie groupoids to the category of Lie algebroids. Let $G$ be a Lie groupoid over $M$ and $g\in G$. The right translation by $g$ is the map $r_g:\source^{-1}(\target(g))\rmap \source^{-1}(\source(g))$ defined by $r_g(h)=hg$. A vector field $X$ on $G$ is said to be \textbf{right invariant} if

$$Tr_gX_h=X_{hg},$$

\noindent for every $(h,g)\in G_{(2)}$. The space $\mathfrak{X}^r(G)$ of right invariant vector fields on $G$ is closed under the bracket of vector fields, so $\mathfrak{X}^r(G)$ becomes a Lie subalgebra of $\mathfrak{X}(G).$

The space of tangent directions to the \textbf{source fibers} is $T^{\source}G=\Ker(T\source)$, which restricts to a vector bundle $AG:=\epsilon^*_MT^{\source}G$ over $M$. Given a section $a$ of $AG$ we define $a^r\in\mathfrak{X}^r(G)$ by $a^r_g=Tr_g(a_{\target(g)})$, yielding an identification $\Gamma_M(AG)\cong\mathfrak{X}^r(G)$. In particular, we have a well defined Lie bracket $[\cdot,\cdot]_{AG}$ on $\Gamma_M(AG)$ given by

\begin{equation}\label{bracketAG}
[a,b]_{AG}^r=[a^r,b^r]
\end{equation}

The vector bundle $AG$ with the bracket \eqref{bracketAG} and the bundle map $\rho_{AG}:=T\target |_{AG}$ make $AG$ into a Lie algebroid over $M$. We say that $AG$ is the Lie algebroid of the Lie groupoid $G$. This construction is functorial in the sense that, given a morphism of Lie groupoids $\Phi:G_1\rmap G_2$ covering a map $\varphi:M_1\rmap M_2$, the tangent map $T\Phi:TG_1\rmap TG_2$ restricts to a bundle map $A(\Phi):AG_1\rmap AG_2$ which defines a Lie algebroid morphism. In this way, we get a functor $A:\mathcal{G}\rmap \mathcal{A}$ from the category of Lie groupoids to the category of Lie algebroids. This functor is referred to as the \textbf{Lie functor}.

We say that a Lie groupoid $G$ is an \textbf{integration} of a Lie algebroid $A$ if $AG$ is isomorphic to $A$, as Lie algebroids. Unlike Lie algebras, the integration of Lie algebroids to Lie groupoids is not possible in general. Obstructions to the existence of integrations were given by Crainic and Fernandes in \cite{CF}. However, if we assume that a Lie algebroid $A$ integrates to a Lie groupoid, then there exists a unique Lie groupoid $G$ with simply connected $\source$-fibers and whose Lie algebroid is $A$. Also Lie's second theorem holds in the category of Lie algebroids.

\begin{theorem}\label{LieII}(Lie's second theorem)

Let $\phi:A_1\rmap A_2$ be a Lie algebroid morphism between integrable Lie algebroids. Let $G_1$ and $G_2$ be Lie groupoids integrating $A_1$ and $A_2$, respectively. If $G_1$ has connected and simply connected $\source$-fibers, then there exists a unique groupoid morphism $\Phi:G_1\rmap G_2$ such that $A(\Phi)=\phi$. 

\end{theorem}

See \cite{CF} for a proof of these facts. 

\subsection{Tangent and cotangent groupoids}\label{tangentcotangentgroupoids}

Given a Lie groupoid $G$ over $M$, we can endow the tangent bundle $TG$ with a Lie groupoid structure over $TM$ as follows. The tangent functor applied to all the structural maps of $G$ gives rise to structural maps $T\source,T\target:TG\rmap TM$, $T\epsilon_M:TM\rmap TG$, $Ti_G:TG\rmap TG$ and $Tm_G:TG_{2}\rmap TG$, satisfying the axioms of a Lie groupoid. We refer to $TG$ with these structural maps as the \textbf{tangent groupoid} associated to $G$. It will be helpful to use the following notation for the tangent multiplication 

$$X_g\bullet Y_h=Tm_G(X_g,Y_h),$$

\noindent for every $(X_g,Y_h)\in (TG)_{(2)}$ composable groupoid pair.

\begin{example}
 
Let $G$ be a Lie group with Lie algebra $\mathfrak{g}$. The tangent bundle $TG$ is also a Lie group, since the space of units consist of only one point. One easily checks that the multiplication on $TG$ is given by

$$X_g\bullet Y_h= T_{g}r_h(X_g) + T_{h}l_g(Y_h).$$

We can use right translations to trivialize $TG$ in such a way that $TG\cong G\times \mathfrak{g}$. With respect to this identification, it is easy to see that the group structure on the tangent bundle corresponds to the semi-direct group $G\ltimes \mathfrak{g}$ determined by the adjoint representation.

\end{example}



Just as the tangent bundle, the cotangent bundle $T^*G$ inherits a Lie groupoid structure over $A^*G$. The source and target maps are defined by

$$\tilde{\source}(\xi_g)a=\xi_g(Tl_g(a-T\target(a)))\quad \text{ and }\quad \tilde{\target}(\eta_g)b=\eta_g(Tr_g(b))$$

\noindent where $\xi_g \in T^*_gG$, $a\in A_{\source(g)}G$ and $\eta_g\in T^*_gG$, $b\in A_{\target(g)}G$. The multiplication on $T^*G$ is defined by

$$(\xi_g\circ \eta_h)(X_g\bullet Y_h)= \xi_g(X_g)+ \eta_h(Y_h)$$

\noindent for $(X_g,Y_h)\in T_{(g,h)}G_{(2)}$.

We refer to $T^*G$ with the groupoid structure over $A^*G$ as the \textbf{cotangent groupoid} of $G$. See \cite{CDW} for more details.

\begin{example}
 
Let $G$ be a Lie group with Lie algebra $\mathfrak{g}$. Then the cotangent groupoid $T^*G$ has base manifold $\mathfrak{g}^*$. We can use right trivializations to identify $T^*G\cong G\times \mathfrak{g}^*$. In terms of this identification, the cotangent groupoid corresponds to the transformation groupoid $G\ltimes\mathfrak{g}^*$ with respect to the coadjoint action.
\end{example}




\begin{remark}
Given a Lie groupoid $G$ over $M$ we have constructed two natural Lie groupoids out of $G$, namely $TG$ over $TM$ and $T^*G$ over $A^*G$. Both the tangent and cotangent Lie groupoids are examples of more general objects called \textbf{$\mathcal{VB}$-groupoids}, that is, Lie groupoids objects in the category of vector bundles. See \cite{Mac-book} for more details.


\end{remark}

\subsection{Tangent and cotangent algebroids}\label{section:tangentcotangentalgebroid}

 Let $A\stackrel{q_A}{\longrightarrow}M$ be a vector bundle. The tangent bundle $TA$ has a natural structure of vector bundle over $TM$, defined by applying the tangent functor to each of the structural maps that define the vector bundle $A\stackrel{q_A}{\longrightarrow}M$. More specifically, we have a square

\begin{align}\label{tangentalgebroid}
\xy 
(-15,10)*+{TA}="t0"; (-15,-10)*+{A}="b0"; 
(15,10)*+{TM}="t1"; (15,-10)*+{M}="b1"; 
{\ar@<.25ex>^{Tq_A} "t0"; "t1"}; 
{\ar@<.25ex>_{q_A} "b0"; "b1"}; 
{\ar@<.25ex>_{p_A} "t0"; "b0"}; 
{\ar@<.25ex>^{p_M} "t1"; "b1"}; 
\endxy
\end{align}

\noindent where all the structural maps of the vector bundle $Tq_A:TA\rmap TM$ are vector bundle morphisms over the corresponding structural maps defining $q_A:A\rmap M$. The vector bundle $Tq_A:TA\rmap TM$ is referred to as the \textbf{tangent prolongation} of $A$. See \cite{Prad, Mac-book} for more details about this construction.

Similarly, the cotangent bundle $T^*A$ determines the following square of vector bundles
 
 \begin{align}\label{cotangentalgebroid}
\xy 
(-15,10)*+{T^*A}="t0"; (-15,-10)*+{A}="b0"; 
(15,10)*+{A^*}="t1"; (15,-10)*+{M}="b1"; 
{\ar@<.25ex>^{r_A} "t0"; "t1"}; 
{\ar@<.25ex>_{q_A} "b0"; "b1"}; 
{\ar@<.25ex>_{c_A} "t0"; "b0"}; 
{\ar@<.25ex>^{q_{A^*}} "t1"; "b1"}; 
\endxy
\end{align}

\noindent where the top arrow $r_A:T^*A\rmap A^*$ is given locally by $r_A(x^i,a^r,p_i,\xi_r)=(x^i,\xi_r)$. Here $(x^i,a^r)\in A_x$, $(x^i,p_i)\in T^*_x M$ and $(x^i,\xi_r)\in A^*_x$, for $x=(x^1,...,x^n)\in M$. The reader can see \cite{Mac-book} for a coordinate free definition. Again, the structural maps of the top horizontal vector bundle are morphisms of vector bundles over the corresponding structural maps of $A$.

\begin{remark}
Both the tangent and cotangent vector bundles are special examples of more general objects called \textbf{double vector bundles}, introduced in \cite{Prad} and further studied by K. Mackenzie \cite{Mac-book}.

\end{remark}

It was shown by Mackenzie and Xu in \cite{Mac-Xu} that whenever $q_A:A\rmap M$ carries a Lie algebroid structure, then $Tq_A:TA\rmap TM$ is a Lie algebroid as well. For that we observe that the space of sections $\Gamma_{TM}(TA)$ is spanned by two types of sections; given a section $a\in\Gamma(A)$ the tangent functor applied to $a$ determines a section $Ta\in\Gamma_{TM}(TA)$. Also $a\in\Gamma(A)$ induces a section $\hat{a}\in\Gamma_{TM}(TA)$ defined by

$$\hat{a}(X)=T0_A(X)+_{p_A} \overline{a(p_M(X))},$$ 

\noindent where $\overline{a(p_M(X))}=\frac{d}{dt}(ta(p_M(X)))\vert_{t=0}$. We refer to $\hat{a}$ as the \textbf{core} section induced by $a$. Now the Lie algebroid structure on $TA$ is defined by

\begin{align*}
[Ta,Tb]_{TA}&=T[a,b]_{A}\\
[Ta,\hat{b}]_{TA}&=\widehat{[a,b]}_A\\
[\hat{a},\hat{b}]_{TA}&=0\\
\rho_{TA}&=J_M\circ T\rho_A
\end{align*}

\noindent Here $J_M:TTM\rmap TTM$ is the canonical involution, which is locally defined by

\begin{equation}\label{canonical-involution}
J_M(x^j,\dot{x}^j, \delta x^j,\delta \dot{x}^j)=(x^j,\delta x^j, \dot{x}^j,\delta \dot{x}^j).
\end{equation}

 \noindent The reader can see \cite{tulc, Mac-book} for more details about this map, including a coordinate free definition.

Similarly, the vector bundle \eqref{cotangentalgebroid} can be endowed with a natural Lie algebroid structure. Concretely, the dual vector bundle $A^*$ inherits a linear Poisson structure. In this case, the cotangent bundle $T^*A^*\rmap A^*$ has a Lie algebroid structure as explained in Example \ref{poissonvsalgebroid}. There is a Legendre type map, which can be defined locally by

\begin{align*}
R_A:T^*&A^*\rmap T^*A\\
(x^i,&\xi_r,p_i,a^r)\mapsto (x^i,a^r,-p_i,\xi_r)
\end{align*}

\noindent One can use this isomorphism to induce a Lie algebroid structure on the vector bundle $r_A:T^*A\rmap A^*$.

We will refer to \eqref{tangentalgebroid} as the \textbf{tangent algebroid} associated to $A$, and \eqref{cotangentalgebroid} will be referred to as the \textbf{cotangent algebroid} of $A$. See \cite{Mac-book, Mac-Xu} for more details about these structures. 

\begin{remark} The Lie algebroids \eqref{tangentalgebroid} and \eqref{cotangentalgebroid} are special examples of double structures. Concretely, each of them  is a \textbf{$\mathcal{LA}$-vector bundle}, that is, a vector bundle object in the category of Lie algebroids. 
\end{remark}

Consider now a Lie groupoid $G$ over $M$. We can apply the Lie functor to the tangent groupoid $TG$ over $TM$, yielding a Lie algebroid $A(TG)\rmap TM$. The canonical involution $J_G:TTG\rmap TTG$ \eqref{canonical-involution} restricts to a Lie algebroid isomorphism

\begin{equation}\label{TA=AT}
 j_G:T(AG)\rmap A(TG). 
 \end{equation}

\noindent Similarly, the Lie functor applied to the cotangent groupoid $T^*G$ gives rise to a Lie algebroid $A(T^*G)\rmap A^*G$. Notice that the canonical pairing $T^*G\times_G TG\rmap \R$ is a groupoid morphism, where $\R$ is viewed as the additive group. Applying the Lie functor to this pairing, we obtain a nondegenerate pairing $A(T^*G)\times_{AG} A(TG)\rmap \R$, as a consequence the vector bundle $A^{\bullet}(TG)$, dual to $A(TG)\rmap AG$, is isomorphic to $A(T^*G)\rmap AG$. The composition of this isomorphism with $j^*_G:A^{\bullet}(TG)\rmap T^*(AG)$, defines a Lie algebroid isomorphism 

\begin{equation}\label{T*A=AT*}
\theta_G:A(T^*G)\rmap T^*(AG).
\end{equation}

\noindent We conclude that, up to canonical identifications, the Lie algebroid of  the tangent groupoid $TG$ (resp. cotangent groupoid $T^*G$) is given by the tangent algebroid $T(AG)$ (resp. cotangent algebroid $T^*(AG)$). The reader can see \cite{Mac-book} for more details about these identifications.


\section{Bivectors and $2$-forms on Lie groupoids}

In this section we study Lie groupoids equipped with bivectors or $2$-forms suitably compatible with the groupoid multiplication. 

\subsection{Multiplicative bivectors}

Let $\pi$ be a Poisson bivector on a smooth manifold $M$ (see Example \ref{poissonvsalgebroid}). We say that a submanifold $Q\subseteq M$ is \textbf{coisotropic} if $\pi^{\sharp}(N^*Q)\subseteq TQ$, where $NQ$ denotes the normal bundle of $Q$.

Consider now a Lie groupoid $G$. Let $\pi_G\in\Gamma(\wedge^2TG)$ be a bivector on $G$. We say that $\pi_G$ is \textbf{multiplicative} if

$$\mathrm{Graph}(m_G)\subseteq G\times G\times \bar{G},$$

\noindent is a coisotropic submanifold, where the last factor on the right hand side is equipped with the Poisson structure $\pi_G$ followed by a minus sign. A \textbf{Poisson groupoid} \cite{weinsteincoisotropic} is a Lie groupoid $G$ equipped with a Poisson bivector $\pi_G$ which is also multiplicative. Poisson groupoids were introduced by A. Weinstein in \cite{weinsteincoisotropic}, generalizing both Poisson Lie groups and symplectic groupoids. This is explained in the following examples.

\begin{example}(Poisson Lie groups)

Let $G$ be a groupoid over a point, i.e. $G$ is a Lie group. A bivector $\pi_G$ on $G$ is multiplicative if and only if the multiplication map $m_G:G\times G\rmap G$ is a Poisson map, which is exactly the definition of a \textbf{Poisson Lie group} \cite{D,LW}.

\end{example}

\begin{example}\label{sympgroupoids}(Nondegenerate case)

Let $(G,\pi_G)$ be a Poisson groupoid. Assume that the bundle map $\pi^{\sharp}_G:T^*G\rmap TG$ is an isomorphism. Then the inverse of $\pi^{\sharp}_G$ defines a $2$-form $\omega_G$ on $G$. Notice that the property of $\pi_G$ being Poisson is equivalent to saying that $\omega_G$ is a closed $2$-form, while the multiplicativity of $\pi_G$ reads

\begin{equation}\label{multiplicativeforms}
m^{*}_G\omega_G=\mathrm{pr}^*_1\omega_G+\mathrm{pr}^*_2\omega_G,
\end{equation}

\noindent where $m_G:G_{(2)}\rmap G$ is the groupoid multiplication and $\mathrm{pr}_1,\mathrm{pr}_2:G_{(2)}\rmap G$ denote the canonical projections. We refer to $(G,\omega_G)$ as a \textbf{symplectic groupoid}.

\end{example}

The reader can see \cite{catfel}, where symplectic groupoids are realized as the phase-spaces of certain sigma models. See also \cite{Wein-Xu} where symplectic groupoids arise in connection with the quantization of Poisson manifolds.

\subsubsection{Infinitesimal description of Poisson groupoids}

\vspace{.2cm}

We would like to express the multiplicativity of a Poisson bivector $\pi_G$ on $G$ in terms of Lie algebroid data. For that one observes that the space of units $M$ of a Poisson groupoid  $G$ defines a coisotropic submanifold of $(G,\pi_G)$. As a consequence, the conormal bundle $N^*M\subseteq T^*G$ is a Lie subalgebroid with respect to the structure defined in \eqref{pisharp} and \eqref{bracketT*M}. Notice that there is a canonical decomposition

$$T_xG=T_xM\oplus A_xG,$$

\noindent for every $x\in M$. As a consequence, the conormal bundle $N^*M$ is canonically isomorphic to $A^*G$, the dual bundle of the Lie algebroid of $G$. Thus, we have a pair of Lie algebroids $(AG,A^*G)$ in duality as vector bundles. These algebroids are compatible in the sense that they define a Lie bialgebroid. 

\begin{definition}\cite{Mac-Xu}

Let $A$ be a Lie algebroid and assume that $A^*$ is also a Lie algebroid. Then the pair $(A,A^*)$ is a Lie bialgebroid if 

$$d_A[\xi,\eta]_{A^*}=[d_A\xi,\eta]_{A^*}+[\xi,d_A\eta]_{A^*},$$

\noindent where $d_A$ is the Lie algebroid differential of $A$.

\end{definition}

Notice that if $A$ is a Lie algebroid over a point, i.e. $A$ is a Lie algebra, then a Lie bialgebroid is nothing but a Lie bialgebra \cite{D}. Just as Lie bialgebras are the infinitesimal data of Poisson Lie groups \cite{D}, Lie bialgebroids may be regarded as the infinitesimal counterpart of Poisson groupoids.\cite{Mac-Xu}.

\begin{remark}\label{pimultiplicative}
It was observed in \cite{Mac-Xu} that the multiplicativity of a bivector $\pi_G$ can be rephrased in terms of Lie groupoid morphisms. In fact, one can check that $\pi_G$ is multiplicative if and only if the bundle map $\pi^{\sharp}_G:T^*G\rmap TG$
is a Lie groupoid morphism covering the anchor $\rho_{A^*G}:A^*G\rmap TM$ of the dual algebroid $A^*G$. Moreover, applying the Lie functor yields a Lie algebroid morphism $A(\pi^{\sharp}_G):A(T^*G)\rmap A(TG)$, which up to the identifications \eqref{TA=AT} and \eqref{T*A=AT*}, coincides with $\pi^{\sharp}_{AG}:T^*(AG)\rmap T(AG)$. Here $\pi_{AG}$ is the linear Poisson bivector on $AG$, induced by the dual Lie algebroid $A^*G$.

\end{remark}

In order to integrate Lie bialgebroids, K. Mackenzie and P. Xu observed in \cite{Mac-Xu} that the Lie bialgebroid condition can be understood in terms of Lie algebroid morphisms. More concretely, the Poisson bivector on the dual of a Lie algebroid $A$ can be described explicitly by considering $(x^1,...,x^m)$ a system of local coordinates on $M$ and $\{e_1,...,e_l\}$ a basis of local sections of $A$, this data determines coordinates $(x^i,a^r)$ on $A$. There are structure functions $\rho^{j}_r$, $C^{t}_{rs}$ for the Lie algebroid $A$, defined by

\begin{itemize}
 \item[i)] $\rho_{A}(e_r)=\rho^{j}_{r}\frac{\partial}{\partial x^j}$, and

 \item[ii)] $[e_r,e_s]_{A}=C^{t}_{rs}e_t.$ 

\end{itemize}

\noindent Now if $\{e^1,...,e^l\}$ is a basis of local sections of $A^*$, dual to $\{e_1,...,e_l\}$, we introduce local coordinates $(x^i,\xi_r)$ on $A^*$. With respect to this local description of $A^*$, the linear Poisson bivector $\pi_{A^*}$ on $A^*$ induced by $A$ has the form

\begin{equation}\label{eq:linearpoissonbivector}
(\pi_{A^*})_{\vert_{(x,\xi)}}=\rho^{i}_{r}(x)\frac{\partial}{\partial x^i}\wedge \frac{\partial}{\partial \xi_r} + \frac{1}{2}C^{t}_{rs}(x)\xi_t\frac{\partial}{\partial \xi_r}\wedge \frac{\partial}{\partial \xi_s}.  
\end{equation}

Here we have used Einstein's convention. The linearity of $\pi_{A^*}$ is equivalent to $\pi^{\sharp}_{A^*}:T^*A^*\rmap TA^*$ being a morphism of double vector bundles \cite{Prad}. It is natural now to find a condition which ensures that $\pi^{\sharp}_{A^*}$ is also a Lie algebroid morphism.

\begin{theorem}\cite{Mac-Xu2}\label{morphicbivector}

If $(A,A^*)$ is a pair of Lie algebroids in duality as vector bundles, then $(A,A^*)$ is a Lie bialgebroid if and only if 

\begin{align}
\xy 
(-15,10)*+{T^*A}="t0"; (-15,-10)*+{A^*}="b0"; 
(15,10)*+{TA}="t1"; (15,-10)*+{TM}="b1"; 
{\ar@<.25ex>^{\pi^{\sharp}_A} "t0"; "t1"}; 
{\ar@<.25ex>_{\rho_{A^*}} "b0"; "b1"}; 
{\ar@<.25ex>_{r_A} "t0"; "b0"}; 
{\ar@<.25ex>^{Tq_A} "t1"; "b1"}; 
\endxy
\end{align}

\noindent is a Lie algebroid morphism between the tangent Lie algebroid \eqref{tangentalgebroid} and the cotangent Lie algebroid \eqref{cotangentalgebroid}. Here the map $\rho_{A^*}$ is the anchor map of the dual Lie algebroid $A^*$. 

\end{theorem}

See \cite{Mac-Xu2} for the proof of this fact. This characterization of Lie bialgebroids combined with remark \ref{pimultiplicative} and Lie's second Theorem \ref{LieII} gives rise to the integration of Lie bialgebroids to Poisson groupoids. See \cite{Mac-Xu, Mac-Xu2} for more details.

\subsection{Multiplicative $2$-forms}

We proceed now to study Lie groupoids $G$ endowed with a closed $2$-form $\omega_G\in\Omega^2(G)$, which is multiplicative as in \eqref{multiplicativeforms}. Notice that Poisson groupoids do not include all multiplicative closed $2$-forms as special cases, since they only include those $2$-forms which are nondegenerate (see example \ref{sympgroupoids}). Most of the results in this subsection can be found in \cite{BCO}.

Notice that a $2$-form on a manifold $M$ determines a skew symmetric bundle map $\omega^{\sharp}:TM\rmap T^*M$ defined by

\begin{equation}\label{omegasharp}
\omega^{\sharp}(X):=\omega(X,\cdot),
\end{equation}

\noindent for every $X\in TM$.

The analogue for $2$-forms of remark \ref{pimultiplicative} reads.

\begin{proposition}\label{omegamultiplicative}\cite{BCO}

Let $G$ be a Lie groupoid over $M$. A $2$-form $\omega_G$ on $G$ is multiplicative if and only if the bundle map $\omega^{\sharp}_G:TG\rmap T^*G$ as in \eqref{omegasharp} is a groupoid morphism between the tangent and cotangent Lie groupoids.
\end{proposition}

If $\omega_G$ is a multiplicative $2$-form on $G$, then the groupoid morphism given by Proposition \ref{omegamultiplicative} covers the bundle map $\lambda:TM\rmap A^*G$, which is the fiberwise dual map of $-\sigma_{\omega_G}:AG\rmap T^*M$, defined by 

\begin{equation}\label{sigma}
\sigma_{\omega_G}(a)=\omega_G(a,\cdot)|_{TM}.
\end{equation}

The application of the Lie functor to the groupoid morphism $\omega^{\sharp}_G:TG\rmap T^*G$ yields a Lie algebroid morphism $A(\omega^{\sharp}_G):A(TG)\rmap A(T^*G)$, which followed by the canonical identifications \eqref{TA=AT} and \eqref{T*A=AT*}, determines a morphism of Lie algebroids

 \begin{align}\label{morphicform}
\xy 
(-15,10)*+{T(AG)}="t0"; (-15,-10)*+{TM}="b0"; 
(15,10)*+{T^*(AG)}="t1"; (15,-10)*+{A^*G}="b1"; 
{\ar@<.25ex>^{\omega^{\sharp}_{AG}} "t0"; "t1"}; 
{\ar@<.25ex>_{-\sigma^{t}_{\omega_G}} "b0"; "b1"}; 
{\ar@<.25ex>_{Tq_{AG}} "t0"; "b0"}; 
{\ar@<.25ex>^{r_{AG}} "t1"; "b1"}; 
\endxy
\end{align}

\noindent where the bottom map is given by the bundle map fiberwise dual to $-\sigma_{\omega_G}:AG\rmap T^*M$. More concretely, the algebroid morphism \eqref{morphicform} is given by 

\begin{equation}\label{pullback-sigma}
\omega^{\sharp}_{AG}=-\sigma^*_{\omega_G}\omega_{can},
\end{equation}

\noindent where $\omega_{can}$ is the canonical symplectic structure on the cotangent bundle $T^*M$ and $\sigma_{\omega_G}$ is given by \eqref{sigma}. A detailed proof of these statements can be found in \cite{BCO}. This discussion motivates the following definition.

\begin{definition}\label{def:morphicform}

A $2$-form $\omega_A$ on a Lie algebroid $A\rmap M$ is called \textbf{morphic} if the induced bundle map $\omega^{\sharp}_A:TA\rmap T^*A$ is a Lie algebroid morphism.

\end{definition}

The $2$-form $\omega_{AG}\in\Omega^2(AG)$ associated to a multiplicative $2$-form $\omega_G$ is morphic. We refer to $\omega_{AG}$ as the \textbf{morphic $2$-form} of $\omega_G$. It turns out that multiplicative $2$-forms on a source simply connected Lie groupoid are in one-to-one correspondence with morphic $2$-forms on its Lie algebroid.

\begin{proposition}\cite{BCO}\label{Lie-multiplicative}

Let $G$ be a Lie groupoid over $M$ with Lie algebroid $A\rmap M$. Assume that $G$ has connected and simply connected source fibers. Then there is a one-to-one correspondence

\begin{align}
\mathrm{Lie}:\Omega^2_{cl, mult}(G)&\rmap \Omega^2_{cl, mor}(A)\\
\omega_G&\mapsto \omega_{AG}
\end{align}

\noindent between multiplicative closed $2$-forms on $G$ and morphic closed $2$-forms on $AG$. Here $\omega_{AG}$ is determined by \eqref{morphicform}. 
\end{proposition}

It is worthwhile to mention here that the correspondence in Proposition \ref{Lie-multiplicative} holds without the closedness assumption. For a detailed discussion of the general case see \cite{BCO}.

Since a multiplicative $2$-form $\omega_G$ induces a map $\sigma_{\omega_G}$ as in \eqref{sigma}, and the morphic $2$-form \eqref{morphicform} associated to $\omega_G$ is related to \eqref{sigma} via \eqref{pullback-sigma}, we would like to understand this relation deeply. For that, let us recall some terminology introduced in \cite{BCWZ}.

\begin{definition}\label{def:IMform}

Let $A$ be a Lie algebroid over $M$. A bundle map $\sigma:A\rmap T^*M$ is called an \textbf{IM-$2$-form} on $A$ if the following identities hold
\begin{align}
\langle \sigma(a), \rho_{A}(b)\rangle=&-\langle \sigma(b),\rho_{A}(a)\rangle\\
\vspace{.2cm}
\sigma[a,b]_{A}=&\Lie_{\rho_{A}(a)}\sigma(b)-\Lie_{\rho_{A}(b)}\sigma(a)+d\langle \sigma(a),\rho_{A}(b)\rangle.
\end{align}

\noindent for every $a,b\in\Gamma(A)$. 

\end{definition}

It was proved in \cite{BCWZ} that the bundle map \eqref{sigma} associated to a \emph{closed} multiplicative $2$-form $\omega_G\in\Omega^2(G)$ defines an IM-$2$-form. This explains why this terminology was introduced, since IM stands for \emph{infinitesimally multiplicative}. It was also shown in \cite{BCWZ} that if $G$ is a Lie groupoid with connected and simply connected source fibers, then there is a one-to-one correspondence between multiplicative closed $2$-forms on $G$ and IM-$2$-forms on the Lie algebroid $AG$ of $G$. The correspondence is given by \eqref{sigma}. The technical point of this results is, as usual, the integration procedure. This is based on the construction of a closed $2$-form on the space of $A$-paths of the Lie algebroid $A$, which is basic with respect to $A$-homotopies, therefore it defines the desired multiplicative closed $2$-form on the $A$-paths model of the Lie groupoid integrating $A$. See Theorem 2.5 of \cite{BCWZ}.

One can avoid $A$-paths and infinite dimensional tools as follows. Given a Lie algebroid $q_A:A\rmap M$, consider the canonical symplectic form $\omega_{can}$ on the cotangent bundle $T^*M$. Every bundle map $\sigma:A\rmap T^*M$ can be used to define a $2$-form on $A$ via

\begin{equation}\label{def:lambda}
\Lambda:=-\sigma^*\omega_{can}.
\end{equation}

This $2$-form is linear in the sense that $\Lambda^{\sharp}:TA\rmap T^*A$ is a morphism of double vector bundles \cite{BCO}. We would like to find a purely infinitesimal condition on $\sigma:A\rmap T^*M$ in such a way that $\Lambda$ be, not only linear, but also a morphic $2$-form. It turns out that the conditions are the ones defining an IM-$2$-form, as established in the following result.

\begin{theorem}\label{integrationsigma}\cite{BCO}

Let $A$ be a Lie algebroid over $M$. Consider a bundle map $\sigma:A\rmap T^*M$ and $\Lambda\in\Omega^2(A)$ given by \eqref{def:lambda}. Then $\sigma$ is an IM-$2$-form if and only if the induced bundle map $\Lambda^{\sharp}:TA\rmap T^*A$ is a Lie algebroid morphism between the tangent Lie algebroid \eqref{tangentalgebroid} and the cotangent Lie algebroid \eqref{cotangentalgebroid}.

\end{theorem}

The morphism $\Lambda^{\sharp}:TA\rmap T^*A$ in Theorem \ref{integrationsigma} covers the bundle map $\lambda:=-\sigma^t:TM\rmap A^*$, which is the fiberwise dual of $\sigma:A\rmap T^*M$. The proof of Theorem \ref{integrationsigma} consists on a long computation involving special sections of the double vector bundles $Tq_A:TA\rmap TM$ and $r_A:T^*A\rmap A^*$. The reader can find a detailed proof in \cite{BCO}.

Combining Lie's second Theorem \ref{LieII} with Proposition \ref{omegamultiplicative} and Theorem \ref{integrationsigma} we recover the correspondence between multiplicative closed $2$-forms and IM-$2$-forms proved previously in \cite{BCWZ}. See \cite{BCO} for this alternative and simpler proof. See also \cite{Bur-Cab} for an extension of these results to multiplicative forms of arbitrary degree.

The following corollary of Theorem \ref{integrationsigma} will be useful later.

\begin{corollary}\label{cor:integrationsigma}
Let $G$ be a Lie groupoid with Lie algebroid $AG$. Assume that $G$ has connected and simply connected source fibers. Then a multiplicative closed $2$-form $\omega_G$ on $G$ is nondegenerate if and only if the associated IM-$2$-form $\sigma_{\omega_G}:AG\rmap T^*M$ is a vector bundle isomorphism.

\end{corollary}

\begin{proof}

It is well known that the nondegeneracy of $\omega_G$ implies that $\sigma_{\omega_G}$ is an isomorphism. See for instance \cite{CDW}. Assume now that $\sigma_{\omega_G}:AG\rmap T^*M$ is an isomorphism. Consider the linear $2$-form on $AG$ defined by $\omega_{AG}=-\sigma^*_{\omega_G}\omega_{can}$. Since $\sigma_{\omega_G}$ is an isomorphism, we conclude that $\omega_{AG}$ is nondegenerate. In particular, the bundle map $\omega^{\sharp}_{AG}:T(AG)\rmap T^*(AG)$ is a vector bundle isomorphism. Since $\sigma_{\omega_G}$ is an IM-$2$-form on $AG$ it follows from Theorem \ref{integrationsigma} that $\omega^{\sharp}_{AG}$ is a Lie algebroid isomorphism which integrates to the Lie groupoid isomorphism $\omega^{\sharp}_{G}$. This shows that $\omega_G$ is nondegenerate.

\end{proof}

\section{Symmetries of Poisson structures}

In this section we study symmetries of Poisson structures. In order to give a conceptually clear exposition we will view Poisson structures as instances of more general geometric structures called Dirac structures.

\subsection{Dirac structures}

Let $M$ be a smooth manifold. Consider the direct sum vector bundle $\mathbb{T}M:=TM\oplus T^*M$, equipped with the nondegenerate symmetric pairing 

$$\langle (X,\alpha),(Y,\beta) \rangle:= \alpha(Y)+\beta(X).$$

\noindent The space of sections $\Gamma(\mathbb{T}M)=\mathfrak{X}(M)\oplus\Omega^1(M)$ is endowed with the Courant bracket

$$\Cour{(X,\alpha),(Y,\beta)}:=([X,Y],\Lie_X\beta -i_Y\mathrm{d}\alpha).$$

A \textbf{Dirac structure} on $M$ \cite{courant} is a sub-bundle $L\subseteq \mathbb{T}M$ which is Lagrangian (i.e. $L^{\perp}=L$) and involutive in the sense that $\Cour{\Gamma(L),\Gamma(L)}\subseteq \Gamma(L)$. Dirac structures were introduced in \cite{courant} motivated by the study of mechanical systems with constraints. The complex analogue yields to the so called generalized complex geometry, see \cite{gualtieri} for further information.

The following examples of Dirac structures will be fundamental along this note.

\begin{example}\label{2formasdirac}

Let $\omega$ be a $2$-form on $M$. One can view $\omega$ as the skew-symmetric bundle map $\omega^{\sharp}:TM\rmap T^*M$ given by \eqref{omegasharp}. The graph $L_{\omega}:=\{(X,\omega^{\sharp}(X))\mid X\in TM\}$ of $\omega^{\sharp}$ is naturally a sub-bundle of $\mathbb{T}M$ which, due to the skew-symmetry of $\omega$, is Lagrangian. The involutivity of $L_{\omega}$ with respect to the Courant bracket is equivalent $\mathrm{d}\omega=0$.

\end{example}

\begin{example}\label{poissonasdirac}

Let $\pi$ be a bivector on $M$ and consider the bundle map $\pi^{\sharp}:T^*M\rmap TM$ defined in \eqref{pisharp}.Then the graph $L_{\pi}:=\{(\pi^{\sharp}(\alpha),\alpha)\mid \alpha\in T^*M\}$ is a Lagrangian sub-bundle of $\mathbb{T}M$. The involutivity condition with respect to the Courant bracket is equivalent to $\pi$ being a Poisson structure, i.e. $[\pi,\pi]=0$, where $[\cdot,\cdot]$ is the Schouten bracket of multivectorfields. 

\end{example}

One observes that a Dirac structure $L$ on $M$ defines a Lie algebroid structure on the vector bundle $L\rmap M$, where the bracket is given by the Courant bracket and the anchor map is the restriction to $L$ of the canonical projection $\mathrm{pr}_{TM}:\mathbb{T}M\rmap TM$. In the case of Dirac structures as in Example \ref{poissonasdirac} the canonical projection $L_{\pi}\rmap T^*M$ is an isomorphism of vector bundles. Since $L_{\pi}$ is a Lie algebroid over $M$, we can use this isomorphism $L_{\pi}\cong T^*M$ to induce a Lie algebroid structure on the cotangent bundle $T^*M$. One easily verifies that this algebroid structure on $T^*M$ coincides with the one defined by \eqref{pisharp} and \eqref{bracketT*M}.

Since Dirac structures define Lie algebroids, we can ask for a Lie groupoid integrating a Dirac structure. For instance, it was proved in \cite{catfel} that if the Lie algebroid $(T^*M)_{\pi}$ associated to a Poisson manifold $(M,\pi)$ is integrable, then the source simply connected Lie groupoid $G$ integrating $(T^*M)_{\pi}$ inherits a symplectic structure $\omega_G$ making the pair $(G,\omega_G)$ into a symplectic groupoid. That is, Poisson manifolds may be regarded as the infinitesimal objects associated with symplectic groupoids. The general problem of the integration of Dirac structures was solved by Bursztyn \emph{et al} in \cite{BCWZ}.

\subsection{Symmetries of $(\mathbb{T}M,\langle\cdot,\cdot\rangle,\Cour{\cdot,\cdot})$}

Given a diffeomorphism $f:M\rmap M$, we can induce a natural automorphism $\mathbb{T}f$ of $\mathbb{T}M$ by $\mathbb{T}f(X,\alpha)=(Tf(X),(Tf^{-1})^*\alpha)$. One easily checks that $\mathbb{T}f$ is an isometry of $(\mathbb{T}M,\langle\cdot,\cdot\rangle)$. A straightforward computation shows that $\mathbb{T}f$ preserves also the Courant bracket. There is another type of symmetries of the data $(\mathbb{T}M,\langle\cdot,\cdot\rangle,\Cour{\cdot,\cdot})$, these are given by the action of closed $2$-forms. More specifically, for any $2$-form $B$ on $M$, define the bundle map

\begin{align}\label{Bfield}
\tau_B:&\mathbb{T}M\rmap \mathbb{T}M\\
& (X,\alpha)\mapsto (X,\alpha+i_X B).
\end{align}

It follows from the skew-symmetry of $B$ that $\tau_B$ preserves the pairing $\langle\cdot,\cdot\rangle$. The following well known result characterizes those $2$-forms $B$ for which $\tau_B$ preserves also the Courant bracket.

\begin{proposition}
 
The bundle automorphism $\tau_B:\mathbb{T}M\rmap \mathbb{T}M$ preserves the Courant bracket if and only if $\mathrm{d}B=0$.

\end{proposition}
 
\begin{proof}

Let $(X,\alpha)$ and $(Y,\beta)$ be sections of $\mathbb{T}M$. Then

$$\Cour{(X, \alpha + i_{X}B),(Y,\beta +i_{Y}B)}=([X,Y], \Lie_{X}\beta - i_Xd\alpha + \Lie_{X}i_{Y}B - i_{Y}di_{X}B),$$

\noindent and using the formula $i_{[X,Y]}=[\Lie_X,i_Y]$ one can see that $B$ is closed if and only if 

$$\Cour{(X, \alpha + i_{X}B),(Y,\beta +i_{Y}B)}=([X,Y], \Lie_{X}\beta - i_Xd\alpha + i_{[X,Y]}B),$$

\noindent which is equivalent to saying that $\tau_{B}$ preserves the Courant bracket.

\end{proof}

Actually, the group of symmetries of $(\mathbb{T}M,\langle\cdot,\cdot\rangle,\Cour{\cdot,\cdot})$ is given by $\mathrm{Diff}(M)\ltimes \Omega_{cl}^2(M)$. See for instance \cite{gualtieri}. As a consequence, the space of Dirac structures on $M$ is invariant under the action of $\mathrm{Diff}(M)\ltimes \Omega_{cl}^2(M)$. If $L$ is a Dirac structure on $M$ and $B\in\Omega^2_{cl}(M)$ then the Dirac structure $\tau_B(L)$ is referred to as a \textbf{gauge transformation} of $L$ by the closed $2$-form $B$. In the literature, closed forms acting as above are also called \textbf{$B$-field transformations}. See for instance \cite{gualtieri, SW}.

\begin{remark}\label{gaugeisomalgebroid}
One immediatly observes that the algebroid structure of a Dirac structure $L$ on $M$ is preserved along its orbit by the action of $\Omega^2_{cl}(M)$.
\end{remark}


\subsection{$B$-field transformations of Poisson manifolds}

Here we are concerned with $B$-fields acting on Poisson manifolds. For that, let $(M,\pi)$ be a Poisson manifold and $B\in\Omega^2(M)$ a closed $2$-form.

\begin{lemma}\label{imformgauge}

The bundle map $\Phi_B:=(\mathrm{Id}+ B^{\sharp}\circ\pi^{\sharp}):(T^*M)_{\pi}\rmap T^*M$ is an IM-$2$-form on the Lie algebroid $(T^*M)_{\pi}$.

\end{lemma}

 \begin{proof}
 
 The tangent bundle $TM$ is a Lie algebroid over $M$, where the anchor map is the identity and the Lie bracket is the usual bracket of vector fields. Notice that a $2$-form $B$ on $M$ is closed if and only if the induced bundle map $B^{\sharp}:TM\rmap T^*M$ is an IM-$2$-form on the Lie algebroid $TM$. Combining this observation with the fact that $\pi^{\sharp}:(T^*M)_{\pi}\rmap TM$ is a Lie algebroid morphism, the result follows.

 \end{proof}

We will see that the IM-$2$-form $\Phi_B$ of Lemma \ref{imformgauge} arises naturally in the study of $B$-field transformations of Poisson structures.

\begin{definition}

Let $(M,\pi)$ be a Poisson manifold. A closed $2$-form $B\in\Omega^2(M)$ is called  \textbf{$\pi$-admissible} if the Dirac structure $\tau_B(L_\pi)$ is also given by the graph of a Poisson bivector $\pi_B\in\Gamma(\wedge^2 TM)$.
\end{definition}

One can easily check that $B$ is $\pi$-admissible if and only if the IM-$2$-form $\Phi_B:(T^*M)_{\pi}\rmap T^*M$ of Lemma \ref{imformgauge} is a vector bundle isomorphism. As a result, the bundle map  $\Phi_B$ is a Lie algebroid isomorphism between $(T^*M)_{\pi}$ and  $(T^*M)_{\pi_B}$. Combining this observation with remark \ref{gaugeisomalgebroid} we conclude that if $(M,\pi)$ integrates to a source simply connected symplectic groupoid $(G,\omega_G)$, then $G$ also integrates the Lie algebroid $(T^*M)_{\pi_B}$. We would like to determine how the symplectic structure on $G$ changes under a $B$-field transformation of the Poisson manifold $(M,\pi)$. The answer is not difficult, since $\Phi_B$ is an IM-$2$-form on the Lie algebroid $(T^*M)_{\pi}$, we can use Theorem \ref{integrationsigma} to integrate $\Phi_B$ to a closed multiplicative $2$-form $\omega^B_G$ on $G$, which is necessarily symplectic due to Corollary \ref{cor:integrationsigma}. In order to give an explicit formula for $\omega^B_G$, we need the following lemma.

\begin{lemma}
Let $G$ be a Lie groupoid over M. For every closed $2$-form $B$ on $M$, the closed $2$-form 
$$B_G:=\target^*B-\source^*B,$$

\noindent is a multiplicative form on $G$ with associated IM-$2$-form $\sigma_{B_G}=B^{\sharp}\circ \rho_{AG}.$

\end{lemma}

\begin{proof}

We need to check that $B_G$ satisfies identity \eqref{multiplicativeforms}. It follows from the definition of $\omega^B_G$ that

\begin{equation}\label{eq:multBG1}
m^*_GB_G=(\target\circ \mathrm{pr}_1)^*B- (\source\circ\mathrm{pr}_2)^*B,
\end{equation}

\noindent where $\mathrm{pr}_1,\mathrm{pr}_2:G_{(2)}\rmap G$ are the canonical projections. One can also see that

\begin{equation}\label{eq:multBG2}
\mathrm{pr}^*_1B_G=(\target\circ\mathrm{pr}_1)^*B- (\source\circ\mathrm{pr}_1)^*B,
\end{equation}

and 

\begin{equation}\label{eq:multBG3}
\mathrm{pr}^*_2B_G=(\target\circ\mathrm{pr}_2)^*B- (\source\circ\mathrm{pr}_2)^*B.
\end{equation}

Since \eqref{multiplicativeforms} is an identity of forms defined on $G_{(2)}$ and $\target\circ\mathrm{pr}_2=\source\circ\mathrm{pr}_1$ on $G_{(2)}$, one combines equations \eqref{eq:multBG1}, \eqref{eq:multBG2} and \eqref{eq:multBG3} to conclude that $B_G$ is multiplicative.

Recall that the IM-$2$-form of a closed multiplicative form is given by \eqref{sigma}. Therefore, the IM-$2$-form of $B_G$ is given at every $a\in AG$ by

\begin{align*}
\sigma_{B_G}(a)=& B_G(a,\cdot)|_{TM}\\
=&B(T\target(a),\cdot)|_{TM} - B(T\source(a),\cdot)|_{TM}\\
=&B(\rho_{AG}(a),\cdot)|_{TM} 
\end{align*}

\noindent where in the last equality we have used that $AG:=\Ker{T\source}|_{M}$.

\end{proof}

The symplectic groupoid $(G,\omega_G)$ defines an IM-$2$-form $\sigma_{\omega_G}:AG\rmap (T^*M)_{\pi}$ which is a Lie algebroid \emph{isomorphism}. For every $a\in AG$ we define $\alpha=\sigma_{\omega_G}(a)$, then

$$\sigma_{B_G}(a)=B^{\sharp}(\rho_{AG}(a))=B^{\sharp}(\pi^{\sharp}(\alpha)).$$

This says that identifying $AG\cong (T^*M)_{\pi}$ via $\sigma_{\omega_G}$, the IM-$2$-form of $B_G=\target^*B-\source^*B$ corresponds to the bundle map $B^{\sharp}\circ\pi^{\sharp}:T^*M\rmap T^*M$. We summarize these observations in the following proposition.

\begin{proposition}

The IM-$2$-form $\Phi_B=(\mathrm{Id}+B^{\sharp}\circ\pi^{\sharp}):(T^*M)_{\pi}\rmap (T^*M)_{\pi_B}$ integrates to the multiplicative closed $2$-form on $G$,  given by

$$\omega^B_G=\omega_G+B_G.$$

Moreover, since $\Phi_B$ is an isomorphism, the pair $(G,\omega^B_G)$ is a symplectic groupoid integrating the Poisson manifold $(M,\pi_B)$.

\end{proposition}

This result was proved first in \cite{BRad} without using IM-$2$-forms. The discussion done in this section provides a conceptually clear proof of some results in \cite{BRad}.



\section{Symmetries of multiplicative Poisson structures}

We want to study $B$-field symmetries of Poisson groupoids. Just as symmetries of Poisson structures are well understood in terms of Dirac geometry, we will see that $B$-field transformations of Poisson groupoids can be studied via Dirac structures suitably compatible with a groupoid multiplication.

\subsection{Multiplicative Dirac structures}

Let $G$ be a Lie groupoid over $M$. One observes that since $TG$ and $T^*G$ are Lie groupoids over $TM$ and $A^*G$, respectively, then the direct sum vector bundle $\mathbb{T}G=TG\oplus T^*G$ inherits a Lie groupoid structure over $TM\oplus A^*G$. The following definition can be found in \cite{Ortiz-Thesis}.

\begin{definition}\label{multiplicativedirac}

A Dirac structure $L_G$ on $G$ is called \textbf{multiplicative} if $L_G$ is a Lie subgroupoid of the direct sum groupoid $\mathbb{T}G$.

\end{definition}

The space of all multiplicative Dirac structures on $G$ will be denoted by $\mathrm{Dir}_{mult}(G)$. Let us see some natural examples.

\begin{example}\label{poissongpds}(Poisson groupoids)

Let $\pi_G$ be a Poisson structure on a Lie groupoid $G$. We have seen that $(G,\pi_G)$ is a Poisson groupoid if and only if $\pi^{\sharp}_G:T^*G\rmap TG$ is a groupoid morphism. Equivalently, the Dirac structure $L_{\pi_G}=\{(\pi^{\sharp}_G(\alpha),\alpha)\mid \alpha\in T^*G\}$ of Example \ref{poissonasdirac} is a subgroupoid. That is, the Dirac structure induced by $\pi_G$ is multiplicative.

\end{example}

\begin{example}\label{presymplecticgpds}(Multiplicative closed $2$-forms)

Let $\omega_G$ be a closed $2$-form on $G$. It follows from Proposition \ref{omegamultiplicative} that $\omega_G$ is multiplicative if and only if $\omega^{\sharp}_G:TG\rmap T^*G$ is a groupoid morphism. Equivalently, the corresponding Dirac structure $L_{\omega_G}=\{(X,\omega^{\sharp}_G(X))\mid X\in TG\}\subseteq \mathbb{T}G$ of Example \ref{2formasdirac} is a subgroupoid. This shows that $L_{\omega_G}$ is a multiplicative Dirac structure.

\end{example}

Multiplicative Dirac structures were introduced by the author in his doctoral thesis \cite{Ortiz-Thesis}. These structures unify both multiplicative closed $2$-forms and Poisson bivectors, offering a conceptually clear framework for studying geometric structures compatible with a group(oid) multiplication. See also \cite{Jotz} where homogeneous spaces for multiplicative Dirac structures are studied.

Let $AG$ be the Lie algebroid of $G$. The tangent and cotangent algebroids induce a Lie algebroid structure on the direct sum $\mathbb{T}(AG)=T(AG)\oplus T^*(AG)$ with base $TM\oplus A^*$. It turns out that the direct sum Lie algebroid $\mathbb{T}(AG)$ is isomorphic to the Lie algebroid of $\mathbb{T}G$, via the map 

$$j^{-1}_G \oplus \theta_G:A(TG)\oplus A(T^*G)\rmap \mathbb{T}(AG),$$

\noindent where $j_G:T(AG)\rmap A(TG)$ and $\theta_G: A(T^*G)\rmap T^*(AG)$ are defined by \eqref{TA=AT} and \eqref{T*A=AT*}, respectively.

\begin{definition}\label{IMdirac}

Let $A$ be a Lie algebroid over $M$. An \textbf{IM-Dirac structure} on $A$ is a Dirac structure $L_A$ which is a Lie subalgebroid of the direct sum Lie algebroid $\mathbb{T}A$.

\end{definition}

The space of all IM-Dirac structures on $A$ will be denoted by $\mathrm{Dir}_{IM}(A)$. The following examples are important.

\begin{example}(Lie bialgebroids)

Let $\pi_A$ be a Poisson bivector on $A$. Then the Dirac structure $L_{\pi_A}$ as in example \ref{poissonasdirac} is an IM-Dirac structure if and only if $\pi^{\sharp}_A:T^*A\rmap TA$ is a Lie algebroid morphism. Due to Theorem \ref{morphicbivector}, an IM-Dirac structure given by the graph of a Poisson bivector is equivalent to $(A,A^*)$ being a Lie bialgebroid.

\end{example}

\begin{example}(IM-$2$-forms)

Let $\omega_A=-\sigma^*\omega_{can}$ be a closed $2$-form on $A$. Then the Dirac structure $L_{\omega_A}$ is an IM-Dirac structure if and only if $\omega_A$ is a morphic $2$-form (see definition \ref{morphicform}).

\end{example}

The following result generalizes the correspondence between Poisson groupoids (resp. multiplicative $2$-forms) and Lie bialgebroids (resp. IM-$2$-forms).

\begin{theorem}\cite{Ortiz-Thesis}\label{multiplicative-IMdirac}

Let $G$ be a Lie groupoid with Lie algebroid $AG$. Assume that $G$ has connected and simply connected source fibers. There is a one-to-one correspondence

\begin{align*}
\mathrm{Lie}:\mathrm{Dir}_{mult}&(G)\rmap \mathrm{Dir}_{IM}(AG)\\
& L_G\mapsto L_{AG}
\end{align*}
\noindent where $L_{AG}:= (j^{-1}_G\oplus \theta_G)(A(L_G))$, and $A(L_G)\subseteq A(\mathbb{T}G)$ is the Lie functor.
\end{theorem}

The reader can find a proof of this result in the author's work \cite{Ortiz-Thesis}.


\subsection{Symmetries of $(\mathbb{TG},\langle\cdot,\cdot\rangle,\Cour{\cdot,\cdot})$}

Let $G$ be a Lie groupoid over $M$ with Lie algebroid $AG$. Throughout this subsection we denote by $B_G$ an \emph{arbitrary} multiplicative closed $2$-form on $G$ with induced morphic $2$-form $B_{AG}$.

\begin{proposition}\label{taumultiplicative}

The bundle map $\tau_{B_G}:\mathbb{T}G\rmap \mathbb{T}G$ defined in \eqref{Bfield} is a groupoid morphism.

\end{proposition}

\begin{proof}

This follows directly from the definition of $\tau_{B_G}$ and Proposition \ref{omegamultiplicative}.

\end{proof}

As a consequence of this proposition we observe that multiplicative closed $2$-forms act on multiplicative Dirac structures via $B$-field transformations.

\begin{corollary}\label{cor:1}

Let $L_G$ be a multiplicative Dirac structure on $G$. Suppose that $B_G$ is a closed multiplicative $2$-form on $G$, then the Dirac structure $L^B_G:=\tau_{B_G}(L_G)$ is also multiplicative.

\end{corollary}

\begin{proof}

This is a straightforward consequence of the proposition above and definition \ref{multiplicativedirac}.

\end{proof}

Now we want to find the infinitesimal data of the multiplicative Dirac structure $\tau_{B_G}(L_G)$ in the sense of Theorem \ref{multiplicative-IMdirac}. For that we need to study $B$-field transformations of the data $(\mathbb{T}A,\langle\cdot,\cdot\rangle)$ for an arbitrary Lie algebroid $A$.  Let $A$ be a Lie algebroid over $M$ and let $B_A$ be an \emph{arbitrary} morphic $2$-form on $A$. 

\begin{proposition}\label{taumorphic}
The bundle map $\tau_{B_A}:\mathbb{T}A\rmap \mathbb{T}A$ is a Lie algebroid morphism.

\end{proposition}

\begin{proof}

This is a straightforward consequence of definition \ref{def:morphicform}.

\end{proof}

\begin{corollary}

Let $L_A$ be an IM-Dirac structure on a Lie algebroid $A$. Suppose that $B_A$ is a morphic closed $2$-form on $A$, then the Dirac structure $L^B_A:=\tau_{B_A}(L_A)$ is also an IM-Dirac structure on $A$.

\end{corollary}

On one hand, as proved in Proposition \ref{Lie-multiplicative}, given a multiplicative closed $2$-form $B_G$ on $G$, we can think of the induced morphic $2$-form $B_{AG}$ on $AG$ as the result of applying the Lie functor to $B_G$. On the other hand, one can apply the Lie functor to the groupoid morphism of Proposition \ref{taumultiplicative}, yielding a Lie algebroid morphism $A(\tau_{B_G}):A(\mathbb{T}G)\rmap A(\mathbb{T}G)$, which composed with the Lie algebroid isomorphism $j^{-1}_G\oplus \theta_G$, gives rise to the Lie algebroid morphism $\tau_{B_{AG}}:\mathbb{T}(AG)\rmap \mathbb{T}(AG)$, where $B_{AG}$ is the morphic $2$-form on $AG$ associated to $B_G$ (See Proposition \ref{Lie-multiplicative}). This gives rise to the main result of this note.

\begin{theorem}\label{thm:1}

Let $G$ be a Lie groupoid over $M$ with Lie algebroid $AG$. Consider a multiplicative closed $2$-form $B_G$ on $G$ with induced morphic $2$-form $B_{AG}$ on $AG$. 
The map $\mathrm{Lie}:\mathrm{Dir}_{mult}(G)\rmap \mathrm{Dir}_{IM}(AG)$ conjugates the $B$-field transformations $\tau_{B_G}$ and $\tau_{B_{AG}}$. That is, 

$$\mathrm{Lie}(\tau_{B_G}(L_G))=\tau_{B_{AG}}(\mathrm{Lie}(L_G)).$$

\end{theorem}

This theorem says that the multiplicative Dirac structure $\tau_{B_G}(L_G)$, obtained by a multiplicative $B$-field transformation of $L_G$, is described infinitesimally by the IM-Dirac structure $\tau_{B_{AG}}(L_{AG})$, defined by the associated morphic $B$-field transformation of $L_{AG}$.


\subsection{$B$-field transformations of Poisson groupoids}

Let $(G,\pi_G)$ be a Poisson groupoid over $M$ with Lie bialgebroid $(AG,A^*G)$. Consider an arbitrary multiplicative $2$-form $B_G$ on $G$ with induced morphic $2$-form $B_{AG}$. Combining Corollary \ref{cor:1} and Theorem \ref{thm:1} we conclude the following.

\begin{proposition}\label{prop:1}
 
The Dirac structure $\tau_{B_G}(L_{\pi_G})$ is a multiplicative Dirac structure, which is infinitesimally described by $\tau_{B_{AG}}(L_{\pi_{AG}})$, where $\pi_{AG}$ is the linear Poisson bivector on $AG$ dual to the Lie algebroid $A^*G$.

\end{proposition}

We will study the special case where $B_G=\source^*B-\target^*B$ for some closed $2$-form $B$ on $M$. 

\begin{remark}\label{induced}

A Lie bialgebroid $(A,A^*)$ over $M$ induces a Poisson structure $\pi$ on $M$ given by $\pi^{\sharp}:=\rho_{A}\circ \rho^*_{A^*}$. This Poisson bivector will be called the Poisson structure \textbf{induced} by $(A,A^*)$. See \cite{Mac-Xu}.

\end{remark}

The Lie bialgebroid $(AG,A^*G)$ determined by $(G,\pi_G)$ induces a Poisson structure $\pi$ on $M$ according to remark \ref{induced}. One can check that a closed $2$-form $B$ on $M$ is $\pi$-admissible if and only if $B_G=\source^*B-\target^*B$ is $\pi_G$-admissible (see \cite{B, BRad}). As a consequence, the Dirac structure $\tau_{B_G}(L_{\pi_G})$ of Proposition \ref{prop:1} is given by the graph of the multiplicative Poisson bivector $\pi^B_G$. We are concerned with the effect of $\tau_{B_G}$ on the Lie bialgebroid of $(G,\pi_G)$.

\begin{lemma}\label{lemma-admissible}

Assume that $G$ has connected and source simply connected source fibers. Then $B_G$ is $\pi_G$-admissible if and only if $B_{AG}$ is $\pi_{AG}$-admissible.

\end{lemma}

\begin{proof}

On one hand, example \eqref{poissongpds} shows that the multiplicativity of $\pi_G$ is equivalent to saying that the map $\pi^{\sharp}_G:T^*G\rmap TG$ is a groupoid morphism between the cotangent and tangent groupoids. Similarly, due to the fact that $B_G$ is multiplicative, one concludes from example \eqref{presymplecticgpds} that $B^{\sharp}_G:TG\rmap T^*G$ is a groupoid morphism. As a result, the map $(\mathrm{Id}+B^{\sharp}_G\circ \pi^{\sharp}_G):T^*G\rmap T^*G$ is a composition of groupoid morphisms. On the other hand, the fact that $B_G$ is $\pi_G$-admissible is equivalent to $(\mathrm{Id}+B^{\sharp}_G\circ \pi^{\sharp}_G):T^*G\rmap T^*G$ be a bijection. Combining these observations we conclude that $(\mathrm{Id}+B^{\sharp}_G\circ \pi^{\sharp}_G):T^*G\rmap T^*G$ is a groupoid isomorphism. The application of the Lie functor, followed by the Lie algebroid isomorphism \eqref{T*A=AT*} gives rise to the Lie algebroid isomorphism $(\mathrm{Id}+B^{\sharp}_{AG}\circ \pi^{\sharp}_{AG}):T^*(AG)\rmap T^*(AG)$. As a consequence, $B_{AG}$ is $\pi_{AG}$-admissible. The converse follows from Lie's second theorem.

\end{proof}

According to Proposition \ref{prop:1}, the infinitesimal data of the multiplicative Poisson bivector $\pi^B_G$ on $G$ corresponds to the Poisson bivector $\pi^{B}_{AG}$ on $AG$, defined by 

$$(\pi^{B}_{AG})^{\sharp}=\pi^{\sharp}_{AG}\circ (\mathrm{Id}+B^{\sharp}_{AG}\circ \pi^{\sharp}_{AG})^{-1}  .$$

Since $(\pi^{B}_{AG})^{\sharp}$ is given by the composition of Lie algebroid morphisms, it follows from Theorem \ref{morphicbivector} that the $B$-field transformation of $\pi_{AG}$ by $B_{AG}$ maps the Lie bialgebroid $(AG,A^*G)$ into another Lie bialgebroid $(AG,(A^*G)_B)$. We will describe the Lie algebroid $(A^*G)_B$ explicitly. For that, we use Theorem \ref{morphicbivector} applied to the Lie algebroid morphism

\begin{align}
\xy 
(-15,10)*+{T^*A}="t0"; (-15,-10)*+{A^*}="b0"; 
(15,10)*+{TA}="t1"; (15,-10)*+{TM}="b1"; 
{\ar@<.25ex>^{(\pi^{B}_{AG})^{\sharp}} "t0"; "t1"}; 
{\ar@<.25ex>_{\rho^{B}_{A^*G}} "b0"; "b1"}; 
{\ar@<.25ex>_{r_A} "t0"; "b0"}; 
{\ar@<.25ex>^{Tq_A} "t1"; "b1"}; 
\endxy
\end{align}

\noindent where $\rho^{B}_{A^*G}:=\rho_{A^*G}\circ (\mathrm{Id}-\rho^*_{AG}\circ B^{\sharp}\circ \rho_{A^*G})^{-1}$, to conclude that the dual vector bundle $A^*G$ inherits a new Lie algebroid structure in such a way that the map $(\mathrm{Id}-\rho^*_{AG}\circ B^{\sharp}\circ \rho_{A^*G}):A^*G\rmap (A^*G)_B$ is a Lie algebroid isomorphism.  This is in agreement with the results of \cite{B}. The Lie bialgebroid $(AG,(A^*G)_B)$ is called the \textbf{gauge transformation} of $(AG,A^*G)$ by the closed $2$-form $B$ on $M$, see \cite{B} for more details.

Using this terminology, the infinitesimal data of the Poisson groupoid $(G,\pi^B_G)$ is given by the gauge transformation of $(AG,A^*G)$ associated to the closed $2$-form $B$ on $M$. Thus we see that Theorem \ref{thm:1} recovers the following result shown in \cite{B}.

\begin{theorem}\label{thm:henrique}\cite{B}

Let $(G,\pi_G)$ be a Poisson groupoid over $M$ with Lie bialgebroid $(AG,A^*G)$ and induced Poisson structure $\pi$ on $M$. Let $B$ be a closed $2$-form on $M$ and let $B_G=\source^*B-\target^*B$ the associated multiplicative $2$-form on $G$. Then $B$ is $\pi$-admissible if and only if $B_G$ is $\pi_G$-admissible. In this case, the Lie bialgebroid of the Poisson groupoid $(G,\pi^B_G)$ is given by $(AG,(A^*G)_B)$.

\end{theorem}

Given a Poisson groupoid $(G,\pi_G)$ over $M$ with induced Poisson structure $\pi$ on $M$ and $B$ a $\pi$-admissible closed $2$-form on $M$, then combining the previous theorem with Lemma \ref{lemma-admissible} one concludes that $B_{AG}$ is $\pi_{AG}$-admissible, where $B_{AG}$ is the morphic form of $B_G=\target^*B-\source^*B$. The gauge transformation Lie bialgebroid $(AG,(A^*G)_B)$ induces a Poisson structure on $M$, given by

$$\pi^{\sharp}_B=\rho_{AG}\circ (\rho^{B}_{A^*G})^*,$$

\noindent where $\rho^{B}_{A^*G}:=\rho_{A^*G}\circ (\mathrm{Id}-\rho^*_{AG}\circ B^{\sharp}\circ \rho_{A^*G})^{-1}$. A straightforward computation shows that $\pi^{\sharp}_B=\tau_B(\pi_M)$. That is, the Poisson structure induced by the gauge transformation Lie bialgebroid $(AG,(A^*G)_B)$ coincides with the gauge transformation, via $B$, of the Poisson structure induced by the Lie bialgebroid $(AG,A^*G)$. 

We have seen that Theorem \ref{thm:henrique} follows from a more general result, namely Theorem \ref{thm:1} where multiplicative Dirac structures played a fundamental role. As we have mentioned, multiplicative Dirac structures unify both multiplicative closed $2$-forms and Poisson bivectors, yielding a natural framework for studying $B$-field transformations of Poisson groupoids. One of the advantages of introducing multiplicative Dirac structures is that they provide a conceptually clear framework, where we can treat all multiplicative structures in a unified manner, often simplifying existing results and proofs.

\begin{remark}

There is another interesting class of multiplicative Dirac structures, namely those given by foliations compatible with a groupoid multiplication, e.g. if $G$ is a Lie groupoid equipped with a free and proper Lie group action by groupoid automorphisms, then the vertical space of this action defines a multiplicative foliation on $G$. It was proved in \cite{ortiz1} that multiplicative foliations on a Lie group $G$ can be explicitly described in terms of connected normal subgroups of $G$. Infinitesimally, if $G$ is a Lie group with Lie algebra $\mathfrak{g}$, then multiplicative foliations on $G$ correspond to Lie subalgebras of $\mathfrak{g}$, which are also ideals. The general case of multiplicative foliations on Lie groupoids is treated in a work in progress with M. Jotz \cite{Jotz-Ortiz}. 
\end{remark}




\begin{thebibliography}{99}

\bibitem{Arias-Crainic} Arias Abad, C., Crainic, M., Representations up to homotopy of Lie algebroids, \textit{Arxiv:0901-0319v1} (2009)



\bibitem{BC}
Bursztyn, H., Crainic, M., Dirac geometry, quasi-Poisson actions and
D/G-valued moment maps.  \textit{J.
Differential Geom. 82, 501-566} (2009)


\bibitem{Bur-Cab} Bursztyn, H., Cabrera, A., Multiplicative forms at the infinitesimal level, \textit{Arxiv:1001-0534} (2010)

\bibitem{BRad} Bursztyn, H., Radko, O., Gauge equivalence of Dirac structures and symplectic groupoids,
\textit{ Ann. Inst. Fourier (Grenoble)  53,   309-337} (2003)


\bibitem{B} Bursztyn, H., On gauge transformations of Poisson structures,
\textit{Quantum Field Theory and Noncommutative Geometry, Lect.
Notes Phys. 662, Springer Verlag 89-112} (2005)


\bibitem{BCO}
Bursztyn, H., Cabrera, A., Ortiz, C., Linear and multiplicative $2$-forms, \textit{Letters in
Mathematical Physics 90, 59-83} (2009)


\bibitem{BCWZ}
Bursztyn, H., Crainic, M., Weinstein, A., Zhu, C., Integration of
twisted Dirac brackets, \textit{ Duke Math. J. 123, 549-607}  (2004)

\bibitem{CW}
Cannas da Silva, A., Weinstein, A., Geometric models for
noncommutative algebras, \textit{Berkeley Mathematics Lecture Notes, 10.
American Mathematical Society, Providence, RI; Berkeley Center for
Pure and Applied Mathematics, Berkeley, CA.} (1999)

\bibitem{catfel}
Cattaneo, A., Felder, G., Poisson sigma models and symplectic
groupoids,  \textit{Quantization of singular symplectic quotients, 198,  61--93,
Progr. Math., Birkh�user, Basel}  (2001)

\bibitem{CDW}
Coste, A., Dazord, P., Weinstein, A., Groupo\"ides
symplectiques, \textit{Publications du D\'epartement de Math\'ematiques.
Nouvelle S\'erie. A, Vol. 2,  i--ii, 1--62, Publ. D\'ep. Math.
Nouvelle S\'er. A, 87-2, Univ. Claude-Bernard, Lyon} (1987)

\bibitem{courant}
Courant, T., Dirac manifolds, \textit{Trans. Amer. Math. Soc. 319, 631-661} (1990)


\bibitem{CF}
Crainic, M., Fernandes, R., Integrability of Lie brackets. \textit{Ann.
of Math. 157, 575-620}   (2003)


\bibitem{D}
Drinfeld, V.,  Hamiltonian structures on Lie groups, Lie
bialgebras and geometric meaning of the classical Yang- Baxter
equations, \textit{Soviet Math. Dokl. 27 (1) 68-71} (1983)


\bibitem{dufour}
Dufour, J.P.,  Zung, N.T., Poisson structures and their normal
forms, \textit{Progress in Mathematics 242, Birkhauser} (2005)


\bibitem{EV}  Etingof, P.,  Varchenko, A., Geometry and classification of solutions of the classical
dynamical Yang-Baxter equation, \textit{Comm. Math. Phys. 192, 77 - 120} (1998)



\bibitem{gualtieri}
Gualtieri, M., Generalized complex geometry, \textit{D. Phil. thesis,
Oxford University} (2003)



\bibitem{Jotz} Jotz, M., A classification theorem for Dirac homogeneous spaces of Dirac Lie groupoids, \textit{Arxiv:1009.0713} (2010)


\bibitem{Jotz-Ortiz}
Jotz, M., Ortiz, C., Multiplicative foliations and their infinitesimal counterparts, \textit{In preparation}  



\bibitem{LW}
Lu, J.H., Weinstein, A., Poisson Lie groups, dressing
transformations, and Bruhat decompositions, \textit{Journal of
Differential Geometry 31, 501-526} (1990)

\bibitem{Mac-book}
Mackenzie, K., General theory of Lie groupoids and Lie algebroids, \textit{Lecture Notes London Math. Soc. 213} (2005)

\bibitem{Mac-Xu}
Mackenzie, K., Xu, P., Lie bialgebroids and Poisson groupoids, \textit{Duke Math. Journal 73, (2) 415-452 } (1994)



\bibitem{Mac-Xu2}
Mackenzie, K., Xu, P., Integration of Lie bialgebroids, \textit{Topology 39, 445-467} (2000)


\bibitem{ortiz1} Ortiz, C., Multiplicative Dirac structures on Lie groups, \textit{C.R. Acad. Sci. Paris Ser. I. 346, (23-24)
1279- 1282} (2008)


\bibitem{Ortiz-Thesis} Ortiz, C., Multiplicative Dirac structures, \textit{Ph.D. Thesis, IMPA, Rio de Janeiro} (2009)


\bibitem{Prad}
Pradines, J., Repr\'esentation des jets non holonomes par des
morphismes vectoriels doubles souds, \textit{C.R. Acad. Sc. Paris, S\'erie A, 278, 1523-1526} (1974)


\bibitem{STS} Semenov-Tian-Shansky, M., Dressing transformation and
Poisson-Lie group actions, \textit{Publ. RIMS, Kyoto University
21, 1237-1260} (1985)


\bibitem{severa}
\v{S}evera, P., Some title containing the words ``homotopy'' and
``symplectic'', e.g. this one.  \textit{Travaux math\'ematiques. Fasc.
XVI, 121--137} (2005)



\bibitem{SW}
\v{S}evera, P., Weinstein, A., Poisson geometry with a $3$-form background, \textit{Prog. Theor. Phys. Suppl. 144, 145-154} (2001)






\bibitem{tulc}
Tulczyjew, W., The Legendre transformation, \textit{Ann. Inst. H.
Poincar\'e, 27, 101-104}


\bibitem{weinsteincoisotropic} Weinstein, A., Coisotropic calculus and Poisson groupoids, \textit{J. Math. Soc. Japan, 40, 705-727} (1988) 

\bibitem{Wein-Xu} Weinstein, A., Xu, P., Extensions of symplectic groupoids and quantization, \textit{J. Reine. Angew. Math. 417, 159-189} (1991)



\end{thebibliography}
\end{document}